\magnification=1200
\centerline{ANALYTICITY OF THE SUSCEPTIBILITY FUNCTION}
\centerline{FOR UNIMODAL MARKOVIAN MAPS OF THE INTERVAL.}
\bigskip\bigskip
\centerline{by Yunping Jiang\footnote{*}{Math. Dept., Queens
College of CUNY, Flushing, NY 11367 and Math. Dept., Graduate
Center of CUNY and Math. Inst., AMSS, CAS. email:
yunqc@forbin.qc.edu} and David Ruelle\footnote{**}{Math. Dept.,
Rutgers University, and IHES. 91440 Bures sur Yvette, France.
email: ruelle@ihes.fr}.}
\bigskip\bigskip\bigskip\bigskip\noindent
    {\leftskip=2cm\rightskip=2cm\sl$\qquad$ Abstract.
    We study the expression (susceptibility)
$$  \Psi(\lambda)
=\sum_{n=0}^\infty\lambda^n\int_I\rho(dx)X(x){d\over dx}A(f^nx)      $$
where $f$ is a unimodal Markovian map of the interval $I$,
and $\rho=\rho_f$ is the corresponding absolutely continuous invariant measure.
We show that $\Psi(\lambda)$ is analytic near $\lambda=1$, where $\Psi(1)$
is formally the derivative of $\int_I\rho(dx)A(x)$ with respect to $f$ in
the direction of the vector field $X$.\par}
\vfill\eject
\null
\bigskip\bigskip
    In a previous note [Ru] the susceptibility function was analyzed for some
    examples of maps of the interval.  The purpose of the present note
    is to give a concise treatment of the general unimodal Markovian case
    (assuming $f$ real analytic).  We hope that it will similarly be possible
    to analyze maps satisfying the Collet-Eckmann condition.
    Eventually, as explained in [Ru], application of a theorem of Whitney [Wh] should
    prove differentiability of the map $f\mapsto\rho_f$ restricted to a suitable set.
\medskip\noindent
{\bf  Setup}
\medskip
Let $I$ be a compact interval of ${\bf R}$ and $f:I\to I$ be real analytic. We assume that there is $c$ in the interior of $I$ such that $f'(0)=0$, $f'(x)>0$ for $x<c$, $f'(x)<0$ for $x>c$, and $f''(c)<0$.  Replacing $I$ by a possibly smaller interval, we assume that $I=[a,b]$ where $a=f^2(c)$, $b=f(c)$.  We assume that the postcritical orbit $P=\{f^nc:n\ge1\}$ is finite: $P=\{p_1,\ldots,p_m\}$; in particular, $f$ is Markovian.  We shall assume that $f$ is {\it analytically expanding} in the sense of Assumption A below; in particular the periodic orbits of $f$ are assumed to be repelling, and therefore $c$ cannot be periodic.  We also assume that $f$ is topologically mixing [this can always be achieved by replacing $I$ by a smaller interval and $f$ by some iterate $f^N$].
\medskip\noindent
{\bf Theorem.}
\medskip
    {\sl Under the above conditions, and Assumption A stated later, there is a unique $f$-invariant probability measure $\rho$ absolutely continuous with respect to Lebesgue on $I$.  If $X$ is real analytic on $I$, and $A\in{\bf C}^1(I)$, then
$$  \Psi(\lambda)
    =\sum_{n=0}^\infty\lambda^n\int_I\rho(dx)X(x){d\over dx}A(f^nx)      $$
extends to a meromorphic function in ${\bf C}$, without pole on $\{\lambda:|\lambda|=1\}$.}
\medskip\noindent
{\bf  Change of variable}
\medskip
The finite set $\{c\}\cup P$ decomposes $I$ into $m$ subintervals $I_j$, with $2m$ endpoints (we ``double'' the endpoints of consecutive subintervals, distinguishing between a $-$ endpoint at the right of an interval, and a $+$ endpoint at the left).  Note that $\eta=\{I_j:j=1,\ldots,m\}$ is a Markov partition for the map $f$.  Consider the critical values of $f^{n}$.  Then for large $n>0$, the set of critical values will be stabilized and is always $P$. We define {\it polar} endpoints as follows:

(1) $p\in P$ is a polar $-$endpoint of an interval in $\eta$ if $p$ is local maximum value of $f^{n}$ for $n$ large.

(2) $v \in P$ is a polar $+$endpoint of an interval in $\eta$ if $p$ is local minimum value of $f^{n}$ for $n$ large.

Every $p\in P$ is a polar $-$ or $+$endpoint and may be both, $c$ is a nonpolar endpoint on both sides.

We define now an increasing continuous map $\varpi: I\to {\bf R}$ so that $J=\varpi I$ is a compact interval. We write $\varpi I_j=J_j$ for $1\leq j\leq m$. Denote by $\omega$ the inverse of $\varpi$. We assume that $\omega|J_j$ extends to a holomorphic function in a complex neighborhood of $J_j$ for $1\leq j\leq m$ and that for $q\in\{c\}\cup P$, $\omega$ has the property
$$  \omega(\varpi q\pm\xi)=\omega(\varpi q)\pm{\xi^2\over2}+O(\xi^4)  $$
if $q$ is a $\pm$polar endpoint, and
$$  \omega(\varpi q\pm\xi)=\omega(\varpi q)\pm\xi+O(\xi^2)  $$
if $q$ is a nonpolar endpoint. [We should really consider disjoint
copies of the $I_j$ and $J_j$, and disjoint neighborhoods of these
in ${\bf C}$ or in a Riemann surface two-sheeted near polar
endpoints. This would lead to notational complications that we
prefer to omit].
\medskip
Applications of this singular change of coordinate have been used
in [Ji1], [BJR], and [Ru]; the reference [Ji2] contains some more
relevant study regarding the method of singular change of
coordinates in one-dimensional dynamical systems. The reader is
encouraged to compare this method with orbifold metrics in [Th,
Chapter 13]. Another relevant application of this method in
complex dynamical systems can be found in [DH].

\medskip
From now on we shall say that $\varpi q$ is a $\pm$polar
(nonpolar) endpoint if $q$ is $\pm$polar (nonpolar).
\medskip\noindent
{\bf The dynamical system viewed after the change of variable.}
\medskip
For any two intervals $I_j, I_k\in \eta$ with $fI_j\supset I_k$,
we define
$$  \psi_{jk}=\varpi\circ(f|I_j)^{-1}\circ(\omega|J_k)  $$
Note that the $\psi_{jk}$ are restrictions of inverse branches of
$g=\varpi\circ f\circ\omega:J\to J$ to intervals in $\eta$. The
function $\psi_{jk}:J_k\to J_j$ extends holomorphically to a
complex neighborhood of $J_k$.  Indeed, note that $(f|I_j)^{-1}$
is holomorphic except if $I_j$ is one of the two intervals around
$c$, in which case the singularity is corrected by $\omega|I_n$,
where $I_n$ is the rightmost interval in $\eta$.  In other cases
$\omega|I_k$ cancels the singularity of $\varpi|I_j$ by our
definition of $\omega$. [Note that $\psi_{jk}'(x) \geq 0$ or $\leq
0$ on $J_{k}$ and may vanish only at an interval endpoint].
\medskip\noindent
{\bf Assumption A.}
\medskip
    {\sl Each $J_k$, for $k=1,\ldots,m$, has a bounded open connected neighborhood $U_k$
    in ${\bf C}$ such that $\psi_{jk}:J_k\to J_j$ extends to a continuous function
    $\psi_{jk}:\bar U_k\to{\bf C}$ holomorphic in $U_k$, and with $\psi_{jk}\bar U_k\subset U_j$.}
\medskip
One checks that the sets $U_k$ can be assumed to be in
$\epsilon$-neighborhoods of the $J_k$.  Also, Assumption A implies
that periodic points for $g$ are strictly repelling.  The
smoothness of $\omega$, $\varpi$ in the interior of subintervals
shows that the same property holds for $f$, apart from interval
endpoints where we however also assume the property to hold:
\medskip
    {\sl The periodic orbits of $f$ are strictly repelling.}
\medskip\noindent
{\bf Markovian graph.}
\medskip
Consider the Markov partition $\eta=\{ I_j\}$. Let us
write $j\succ k$ ($j$ covers $k$) if $fI_j\supset I_k$ (we allow $j\succ j$).  This defines a directed graph with vertex set $\{1,\ldots,m\}$ and oriented edges $(j,k)$ for $j\succ k$.  Since we have assumed our dynamical system $f$ to be topological mixing, our graph is also mixing in the sense that there is $N\ge1$ such that for all $j,k\in\{1,\ldots,m\}$ we have $j\succ\ldots\succ k$ ($N$ edges) corresponding to $f^NI_j\supset I_k$.

\medskip\noindent
{\bf Transfer operators.}
\medskip
    For a function $\Phi=(\Phi_j)$ defined on $\sqcup J_j$, write
$$
({\cal L}\Phi)_k(z)
=\sum_{j:j\succ k}{\rm sgn}(j)\psi'_{jk}(z)\Phi(\psi_{jk}z)
$$
$$
({\cal L}_0\Phi)_k(z)
=\sum_{j:j\succ k}{\rm sgn}(j)\Phi(\psi_{jk}z)
$$
where ${\rm sgn}(j)$ is $+1$ if $\psi_{jk}$ is increasing on $J_k$, and $-1$ if $\psi_{jk}$ is decreasing on $J_k$. If $H$ is the Hilbert space of functions on
$\sqcup_{j\in L}\bar U_j$ which are square integrable with respect
to Lebesgue, and have holomorphic restrictions to the $U_j$, then
${\cal L}$ and ${\cal L}_0$ acting on $H$ are holomorphy
improving, hence compact and of trace class.
\medskip\noindent
{\bf Properties of ${\cal L}$.}
\medskip
    For $x\in J_k$ we have
$$  ({\cal L}\Phi)_k(x)=\sum_{j\succ k}|\psi'_{jk}(x)|\Phi_j(\psi_{jk}x)  $$
hence $\Phi\ge0$ implies ${\cal L}\Phi\ge0$ (${\cal L}$ preserves positivity) and
$$  \int_Jdx\,({\cal L}\Phi)(x)=\sum_k\int_{J_k} dx({\cal L}\Phi)_k(x)
    =\sum_j\int_{J_j} dx\Phi_j(x)=\int_Jdx\,\Phi(x)      $$
(${\cal L}$ preserves mass).  Using mixing one finds that ${\cal
L}$ has a simple eigenvalue $\mu_0=1$ corresponding to an
eigenfunction $\sigma_0>0$. The other eigenvalues $\mu_\ell$
satisfy $|\mu_\ell|<1$, and their (generalized) eigenfunctions
$\sigma_\ell$ satisfy $\int_Jdx\,\sigma_\ell(x)=0$. If we
normalize $\sigma_0$ by $\int_Jdx\,\sigma_0(x)=1$, then
$\sigma_0(dx)=\sigma_0(x)dx$ is the unique $g$-invariant
probability measure absolutely continuous with respect to Lebesgue
on $J$. In particular, $\sigma_0(x)dx$ is ergodic.
\medskip
Let now $H_1\subset H$ consist of those $\Phi=(\Phi_k)$ such that
the derivative $\Phi'$ vanishes at the (polar) endpoints $\varpi a$ ,$\varpi b$ of $J$, and such that at the common endpoint
$\varpi q$ ($q\in\{c\}\cup P\backslash\{a,b\}$)) of two subintervals we have equality on both
sides of a quantity which is either

$\bullet$ the value of $\Phi$ for a nonpolar endpoint, or

$\bullet$ the value of $\pm\Phi'$ for a polar $\pm$ endpoint.
\medskip
We note that ${\cal L}H_1\subset H_1$ [this requires a case by
case discussion].  Furthermore $\sigma_0\in H_1$ [take $\phi\in H$
such that $\phi\ge0$, $\int_Jdy\,\phi(y)=1$, and $\phi,\phi'$
vanish at subinterval endpoints; then $\phi\in H_1$ and
$\sigma_0=\lim_{n\to\infty}{\cal L}^n\phi\in H_1$].
\medskip\noindent
{\bf Evaluating $\Psi(\lambda)$.}
\medskip
The image $\rho(dx)=\rho(x)dx$ of $\sigma_0(y)dy$ by $\omega$ is
the unique $f$-invariant probability measure absolutely continuous
with respect to Lebesgue on $I$.  We have
$$  \rho(x)=\sigma_0(\varpi x)\varpi'(x)      $$
Consider now the expression
$$  \Psi(\lambda)
    =\sum_{n=0}^\infty\lambda^n\int_I\rho(dx)X(x){d\over dx}A(f^nx)      $$
where we assume that $X$ extends to a holomorphic function in a neighborhood
of each $I_k$ and takes the same value at both sides of common endpoints of intervals in $\eta$ (continuity).
Also assume that $A\in{\cal C}^1(I)$.  For sufficiently small $|\lambda|$,
the series defining $\Psi(\lambda)$ converges.  Writing $B=A\circ\omega$ ($B$ has piecewise continuous derivative) and $x=\omega y$ we have
$$  X(x){d\over dx}A(f^nx)=X(\omega y){1\over\omega'(y)}{d\over dy}B(g^ny)  $$
hence
$$  \Psi(\lambda)=\sum_{n=0}^\infty\lambda^n
\int_Jdy\,\sigma_0(y){X(\omega y)\over\omega'(y)}{d\over dy}B(g^ny)      $$
Defining $Y(y)=\sigma_0(y)X(\omega y)/\omega'(y)$, we see that $Y$ extends
to a holomorphic function in a neighborhood of each $J_k$, which we may take to be $U_k$,
except for a simple pole at each polar endpoint of $J_k$.
Since $\sigma_0\in H_1$, the properties assumed for $\omega$ imply
that also $(X\circ\omega)\times\sigma_0\in H_1$.
Note that near a nonpolar subinterval endpoint $\varpi q$
$$  \omega'(\varpi q\pm\xi)=1+O(\xi)      $$
and near a $\pm$ polar endpoint
$$  \omega'(\varpi q\pm\xi)=\xi+O(\xi^3)      $$
Therefore
$$  Y(\varpi q\pm\xi)=A^{\pm}{1\over\xi}+B^{\pm}+O(\xi)      $$
where $B^+=B^-$ for the two sides of $\varpi q$, and $B^+=0$ at the left endpoint $\varpi a$ of $J$, $B^-=0$
at the right endpoint $\varpi b$ of $J$.  We may write
$$  \int_Jdy\,\sigma_0(y){X(\omega y)\over\omega'(y)}{d\over dy}B(g^ny)
    =\int_Jdy\,Y(y)g'(y)\cdots g'(g^{n-1}y)B'(g^ny)      $$
$$  =\int_Jds\,({\cal L}_0^nY)(s)B'(s)      $$
where ${\cal L}_0$ has been defined above, and we have thus
$$  \Psi(\lambda)
    =\sum_{n=0}^\infty\lambda^n\int_Jds\,({\cal L}_0^nY)(s)B'(s)      $$
\noindent
{\bf Properties of ${\cal L}_0$.}
\medskip
We let now $H_0\subset H$ be the space of functions
$\Phi=(\Phi_k)$ vanishing at the endpoints $\varpi a$, $\varpi
b$ of $J$, and such that the values of $\Phi$ on both sides of
common endpoints of intervals $J_j$ coincide (continuity).
Therefore ${\cal L}_0H_0\subset H_0$.
\medskip
There is a periodic orbit $\gamma_1,\ldots,\gamma_p$ (with
$g\gamma_j=\gamma_{j+1({\rm mod}p)}$) of polar endpoints where
$\gamma_\alpha$ is the $\pm$ endpoint of some subinterval
$J_{k(\alpha)}$.  Choose $P_\alpha$ to be $0$ on subintervals
different from $J_{k(\alpha)}$, and to be holomorphic on a complex
neighborhood of $J_{k(\alpha)}$ except at $\gamma_\alpha$.  Also
assume that
$$  P_\alpha(\gamma_\alpha\pm\xi)={1\over\xi}+O(\xi)      $$
and that $P_\alpha$ vanishes at the endpoint of $J_{k(\alpha)}$ different from $\gamma_\alpha$.  Then
$$  {\cal L}_0P_\alpha-|f'(\gamma(\alpha)|^{1/2}P_{\alpha+1({\rm mod}p)}
    \in H_0      $$
Therefore ${\cal L}_0^pP_1-\Lambda P_1=u\in H_0$ where $\Lambda=\prod_{\alpha=1}^p|f'(\gamma(\alpha)|^{1/2}>1$.  Since the spectrum of ${\cal L}$ acting on $H$ is contained in the closed unit disk, and since the derivative $u'$ is in $H$, we may define $v=({\cal L}^p-\Lambda)^{-1}u'\in H$.  Since $\int_Jdy\,u'(y)=0$ we also have $\int_Jdy\,v(y)=0$ and we can take $w\in H_0$ such that $w'=v$.  We have thus
$$  (({\cal L}_0^p-\Lambda)w)'=({\cal L}^p-\Lambda)w'
    =({\cal L}^p-\Lambda)v=u'      $$
so that $({\cal L}_0^p-\Lambda)w=u$ [there is no additive constant of integration since $({\cal L}_0^p-\Lambda)w$ and $u$ are in $H_0$].  Finally
$$  ({\cal L}_0^p-\Lambda)(P_1-w)=0      $$
There is thus a ${\cal L}_0$-invariant $p$-dimensional vector space spanned by vectors $P_\alpha-w_\alpha$ with $w_\alpha\in H_0$, such that the spectrum of ${\cal L}_0$ restricted to this space consists of eigenvalues $\omega_\ell$ with
$$  \omega_\ell=\Lambda^{1/p}e^{2\pi\ell i/p}
=|\prod_{\alpha=1}^pf'(\gamma_\alpha)|^{1/2p}e^{2\pi\ell i/p}      $$
for $\ell=0,\ldots,p-1$.
\medskip
    For the postcritical but nonperiodic polar points $\tilde\gamma_1,\ldots,\tilde\gamma_q$ define $\tilde P_\beta$ like $P_\alpha$ above, with $\gamma_\alpha$ replaced by $\tilde\gamma_\beta$.  For each $\beta$ there is $\alpha=\alpha(\beta)$ with
$$  {\cal L}_0^q(\tilde P_\beta-\tilde\Lambda_\beta P_\alpha)\in H_0      $$
with some $\tilde\Lambda_\beta\ne0$, hence
$$  {\cal L}_0^q(\tilde P_\beta-\tilde\Lambda_\beta(P_\alpha-w_\alpha))
    =\tilde Y_\beta\in H_0      $$
\noindent
{\bf Poles of $\Psi(\lambda)$.}
\medskip
    We may now write
$$  Y=Y_0+Y_1+Y_2      $$
where
$$  Y_0\in H_0      $$
$$  Y_1=\sum_{\alpha=1}^pc_\alpha(P_\alpha-w_\alpha)      $$
$$  Y_2=\sum_{\beta=1}^q\tilde c_\beta
(\tilde P_\beta-\tilde\Lambda_\beta(P_{\alpha(\beta)}-w_{\alpha(\beta)}))      $$
and there is a corresponding decomposition $\Psi(\lambda)=\Psi_0(\lambda)+\Psi_1(\lambda)+\Psi_2(\lambda)$.  Here $\Psi_1(\lambda)$ is a sum of terms $C_\ell/(\lambda-\omega_\ell)$ where $\omega_\ell=\Lambda^{1/p}\times p$-th root of unity; $\Psi_2(\lambda)=$ polynomial of degree $q-1$ in $\lambda$ plus $\lambda^q\sum_{\beta=1}^q\tilde c_\beta\tilde\Psi_\beta(\lambda)$ where $\tilde\Psi_\beta$ is obtained if we replace $Y$ by $\tilde Y_\beta$ in the definition of $\Psi$.  The poles of $\Psi(\lambda)$ are thus those of $\Psi_1(\lambda)$ at the $\omega_\ell$ and those of $\Psi_0(\lambda)$ and $\tilde\Psi_\beta(\lambda)$.  The discussion is the same for $\Psi_0$ and the $\tilde\Psi_\beta$, we shall thus only consider $\Psi_0$.  Since $Y_0\in H_0$ and ${\cal L}_0H_0\subset H_0$ we have
$$  \Psi_0(\lambda)
    =\sum_{n=0}^\infty\lambda^n\int_Jds\,({\cal L}_0^nY_0)(s)B'(s)
    =-\sum_{n=0}^\infty\lambda^n\int_Jds\,({\cal L}_0^nY_0)'(s)B(s)      $$
$$  =-\sum_{n=0}^\infty\lambda^n\int_Jds\,({\cal L}^nY'_0)(s)B(s)      $$
It follows that $\Psi_0(\lambda)$ extends meromorphically to ${\bf C}$ with poles at the $\mu_k^{-1}$.  We want to show that the residue of the pole at $\mu_0^{-1}=1$ vanishes.  Since $\int_Jdy\,\sigma_k(y)=0$ for $k\ge1$, the coefficient of $\sigma_0$ in the expansion of $Y'_0$ is proportional to
$$  \int_Jdy\,Y'_0(y)=Y(\varpi b)-Y(\varpi a)=0      $$
because $Y_0\in H_0$.  Therefore $\Psi_0(\lambda)$ is holomorphic for $|\lambda|=1$, and the same holds for the $\tilde\Psi_\beta(\lambda)$, concluding the proof of the theorem.  In fact we know that the poles of $\Psi(\lambda)$ are located at $\mu_k^{-1}$ for $k\ge1$, and at $\omega_\ell^{-1}$ for $\ell=0,\ldots,p-1$, so that $|\mu_k^{-1}|>1,|\omega_\ell^{-1}|<1$.
\medskip\noindent
{\bf Acknowledgments.}
\medskip
    One of us (Y.J.) is partially supported by grants from the NSF and the PSC-CUNY Program, and by the Hundred Talents Program of the CAS.
\vfill\eject
{\bf References.}

[BJR] V. Baladi, Y. Jiang, H.H. Rugh.  ``Dynamical determinants
via dynamical conjugacies for postcritically finite polynomials.''
J. Statist. Phys. {\bf 108}, 973-993(2002).

[DH] D. Douady and J. H. Hubbard, ``A proof of Thurston's
topological characterization of rational functions.'' Acta Math.,
Vol. {\bf 171}, 263-297(1993).

[Ji1] Y. Jiang, ``On Ulam-von Neumann transformations.'' Commun.
Math. Phys., {\bf 172}, 449-459(1995).

[Ji2] Y. Jiang, ``Renormalization and Geometry in One-Dimensional
and Complex Dynamics.'' Advanced Series in Nonlinear Dynamics,
Vol. 10 (1996), World Scientific Publishing Co. Pte. Ltd., River
Edge, NJ.

[Ru] D. Ruelle.  ``Differentiating the absolutely continuous
invariant measure of an interval map $f$ with respect to $f$.''
Commun. Math. Phys. to appear.

[Th] W. Thurston, ``The Geometry and Topology of
Three-Manifolds'', Electronic version 1.1-March 2002,
http://www.msri.org/publications/books/gt3m/

[Wh] H. Whitney.  ``Analytic expansions of differentiable
functions defined in closed sets.''  Trans. Amer. Math. Soc. {\bf
36},63-89(1934).

\end